\documentclass[11pt]{article}
\usepackage{amsfonts}
\usepackage{color}
\usepackage{graphics}

\usepackage{cite}
\usepackage{latexsym}
\usepackage{amsmath}
\usepackage{amssymb}
\usepackage[dvips]{epsfig}
\usepackage{amscd}

 \usepackage{color}
 \usepackage{indentfirst}
\usepackage{ntheorem}
\theoremstyle{plain}
\theoremseparator{.}
\newtheorem{theorem}{Theorem}[section]
\newtheorem{remark}{Remark}[section]

\newtheorem{lemma}[theorem]{Lemma}

\newcommand\thmref[1]{Theorem~\ref{#1}}
\newcommand\lemref[1]{Lemma~\ref{#1}}

\newcommand{\nb}{{\bold b}}
\newcommand{\nw}{{\bold w}}
\newcommand{\ti}{\tilde}

\newcommand{\lm}{\lambda}

\def\pf{{\it Proof.}  }

\newcommand{\thatsall}{\hfill$\Box$}
\newcommand{\bi}{\bibitem}

\newcommand{\bt}{\begin{theorem}}
\newcommand{\bl}{\begin{lemma}}
\newcommand{\el}{\end{lemma}}
\newcommand{\et}{\end{theorem}}

\renewcommand{\b}{\beta  }

\newcommand{\te}{\theta}

\newcommand{\al}{\alpha}

\newcommand{\ve}{\varepsilon}
\newcommand{\la}{\label}
\newcommand{\si}{\sigma}
\newcommand{\ka}{\kappa}

\newcommand{\bn}{\begin{eqnarray}}
\newcommand{\en}{\end{eqnarray}}
\newcommand{\bnn}{\begin{eqnarray*}}
\newcommand{\enn}{\end{eqnarray*}}

\newcommand{\bnnn}{\begin{eqnarray*}}
\newcommand{\ennn}{\end{eqnarray*}}
\newcommand{\ben}{\begin{enumerate}}
\newcommand{\een}{\end{enumerate}}
\newcommand{\ba}{\begin{aligned}}
\newcommand{\ea}{\end{aligned}}
\newcommand{\be}{\begin{equation}}
\newcommand{\ee}{\end{equation}}

\def\norm[#1]#2{\|#2\|_{#1}}

\def\xix{\int_0^1}

\makeatletter
\@addtoreset{equation}{section}
\makeatother

\makeatletter
\@addtoreset{equation}{section}
\makeatother

\title{ Global Strong Solutions to   Magnetohydrodynamics  with Density-Dependent Viscosity
and Degenerate Heat-Conductivity  \thanks{ Partially supported by NNSFC   11671027 and 11471321.}
}

\author{Bin Huang,  Xiaoding Shi, Ying Sun \thanks{   Email addresses: abinhuang@gmail.com, abinhuang36@163.com (B. Huang), shixd@mail.buct.edu.cn (X. Shi), 1913349041@qq.com (Y. Sun)}  \\[3mm]     Department of Mathematics, Faculty  of Science, \\Beijing University of Chemical Technology, \\ Beijing  100029, P. R. China }
\date{ }

\begin{document}
\maketitle
\begin{abstract}
We deal with the equations of a planar magnetohydrodynamic compressible flow  with  the viscosity depending on the specific volume of the gas  and the heat conductivity    proportional to a positive power of  the temperature. Under the same   conditions   on  the  initial data  as those of the constant viscosity and   heat conductivity case ([Kazhikhov
(1987)], we obtain the global  existence and uniqueness of strong solutions which means   no shock wave, vacuum, or mass or heat concentration will be developed in finite time, although the motion of the flow has large oscillations and the  interaction between the hydrodynamic and magnetodynamic effects is complex.
Our result can be regarded as a natural generalization of the Kazhikhov's theory  for the constant viscosity and heat conductivity case to that of nonlinear viscosity and  degenerate    heat-conductivity.  \end{abstract}

$\mathbf{Keywords.}$ Magnetohydrodynamics,  Large initial data, Global strong solutions, Degenerate heat-conductivity,  Density-dependent viscosity

{\bf Math Subject Classification:} 35Q35; 76N10.

\section{Introduction}

Magnetohydrodynamics (MHD), concerning the motion of conducting fluids in an electromagnetic field,
covers a wide range of physical objects from liquid metals to cosmic plasmas (\cite{q1,q2,q3,q4,q5,q6,q7}). The dynamic motion of the fluids and the magnetic field interact strongly with each other. Moreover, the hydrodynamic and electrodynamic effects are coupled.
 We are concerned with the governing equations of a planar magnetohydrodynamic compressible flow written
  in the Lagrange variables
  \be\la{1.1}
v_t=u_{x},
\ee
\be\la{1.2}
u_{t}+(P+\frac12 |\nb|^2)_{x}=\left(\mu\frac{u_{x}}{v}\right)_{x},\ee
 \be\la{1.3}\nw_t-\nb_x=\left(\lm\frac{\nw_{x}}{v}\right)_{x}, \ee \be\la{1.4} (v\nb)_t-\nw_x =\left(\nu\frac{\nb_{x}}{v}\right)_{x}, \ee
\be\la{1.5}\ba&\left(e+\frac{u^2+|\nw|^2+v|\nb|^2}{2}\right)_{t}+ \left(
u\left(P+\frac12|\nb|^2\right)-\nw\cdot\nb\right)_{x}\\&\quad=\left(\kappa\frac{\theta_{x}}{v}
+\mu\frac{uu_x}{v}+\lm\frac{\nw\cdot\nw_x}{v}+\nu\frac{\nb\cdot\nb_x}{v}\right)_x,
 \ea
\ee
where $t>0$ is time, $x\in \Omega=(0,1)$ denotes the
Lagrange mass coordinate,  and the unknown functions $v>0, u, \nw\in \mathbb{R}^2, \nb\in   \mathbb{R}^2, e>0, \theta>0$ and $P$ are,  respectively, the specific volume of the gas, longitudinal velocity, transverse velocity, transverse magnetic field, internal energy,  absolute temperature and  pressure. $\mu$ and $\lm$ are the viscosity of the flow,
$\nu$ is the magnetic diffusivity of the magnetic field, and $\ka$ is the heat conductivity.

In this paper, we
concentrate on  a perfect gas for magnetohydrodynamic flow, that is, $P$ and $e$ satisfy \be \la{1.6}   P =R \theta/{v},\quad e=c_v\theta +\mbox{const},
\ee
where  both specific gas constant  $R$ and   heat
capacity at constant volume $c_v $ are   positive constants.
We also assume that $\lm $ and $\nu$ are positive constants, and  $\mu, \ka$ satisfy \be\la{1.7}   \mu=\ti\mu_1+\ti\mu_2 v^{-\al} , \quad \ka=\ti\ka \te^\beta, \ee
with constants $\ti\mu_1>0,\ti \mu_2\ge 0,\ti\ka>0,$ and $\al,\beta\ge 0.$

The system \eqref{1.1}-\eqref{1.7} is supplemented
with the initial  conditions
\be\la{1.8}(v,u,\te,\nb,\nw)(x,0)=(v_0,u_0,\te_0,\nb_0,\nw_0)(x),  \quad  x\in\Omega, \ee
and boundary  conditions\be\la{1.9} \left(u,\nb,\nw,\theta_x\right)|_{\partial\Omega}=0,\ee
where  the initial data \eqref{1.8} should be compatible with the boundary conditions \eqref{1.9}.

There is huge literature on the studies of the global existence and large time behavior of solutions to the  compressible Navier-Stokes system and MHD. Indeed,
for compressible Navier-Stokes system,  Kazhikhov and Shelukhin
\cite{9} first obtained   the global existence of solutions
 for constant coefficients $(\al=\beta=0)$  with large initial data. From then on, much effort has been made to generalize this approach to other cases. Jenssen-Karper \cite{24}  proved the global existence of   weak solutions
under the assumption that $\al=0$ and $\beta\in(0,3/2).$ Later, for $\al=0$ and $\beta\in(0,\infty),$  Pan-Zhang \cite{28} obtained the global strong solutions under the condition that \bnn( v_0,u_0,\te_0)\in H^1\times H^2\times H^2, \enn which was further relaxed to \bnn( v_0,u_0,\te_0)\in H^1, \enn by Huang-Shi \cite{hs1} where they also obtained the large-time behavior of the strong solutions. As for MHD, the existence and uniqueness of
local smooth solutions was first proved in \cite{vh}.
 For constant coefficients $(\al=\beta=0)$  with large initial data,  Kazhikhov  \cite{ka1} (see also \cite{az1})
 first obtained   the global existence of strong solutions.    From then on,
significant progress has been made on the mathematical aspect of the initial and
initial
boundary value problems, see \cite{6,cw1,cw2,fjn1,25,fhl1} and the references therein.      However, it should be mentioned here that the methods used there rely heavily on the non-degeneracy of both the viscosity $\mu$ and the heat conductivity $\ka$  and  cannot be applied directly to  the degenerate and nonlinear case  ($\al\geq 0,\b>0$).

More recently,    Hu-Ju \cite{7}  extended  Pan-Zhang's result (\cite{28}) to the MHD case and   proved the  global strong solutions to the initial-boundary-value problem \eqref{1.1}-\eqref{1.9} with $\al=0$ and $\b>0$ under the condition that \be \la{h1.8}v_0\in H^1,\quad ( u_0,\te_0,\nb_0,\nw_0)\in H^2, \ee which is stronger than that of Kazhikhov (\cite{ka1}). In fact,
 the main aim of this paper is to generalize  Kazhikhov's result \cite{ka1} to  the degenerate and nonlinear case  and  prove the global existence of strong solutions to \eqref{1.1}-\eqref{1.9} with $\al\ge 0,\b>0 $  and \bnn(v_0, u_0,\te_0,\nb_0,\nw_0)\in H^1. \enn

   Then  we state
  our main result   as follows.
 \begin{theorem}\la{thm1.1} Suppose that \be \la{a1.10}\al\ge 0,\quad \beta> 0,\ee
 and that the initial data $ ( v_0,u_0,\te_0,\nb_0,\nw_0)$   satisfies
  \be  \la{a1.11} ( v_0,\te_0)\in   H^1 (0,1),\quad (u_0,\nb_0,\nw_0)\in H^1_0 (0,1),\ee  and \be \la{a1.12}
\inf_{x\in (0,1)}v_0(x)>0, \quad \inf_{x\in (0,1)}\theta_0(x)>0. \ee
Then, the initial-boundary-value problem \eqref{1.1}-\eqref{1.9} has a unique strong solution $(v,u,\te,$\\$\nb,\nw)$ such that for each fixed $ T>0 $,
 \be
 \begin{cases}   v ,\,\theta \in L^\infty(0,T;H^1(0,1)),\quad u,\,\nb,\nw  \in L^\infty(0,T;H^1_0(0,1)),\\ v_t\in
  L^\infty(0,T;L^2(0,1))\cap L^2(0,T;H^1(0,1)), \\ u_t,\,\theta_t,\,\nb_t,\,\nw_t,\,u_{xx},\,\te_{xx},\,\nb_{xx},\,\nw_{xx} \in
  L^2((0,1)\times(0,T)),\end{cases}\ee
  and for each $(x,t)\in[0,1]\times[0,T]$
  \be C^{-1}\leq v(x,t)\leq C,\quad C^{-1}\leq\te(x,t)\leq C,\ee
  where $C>0$ is a constant  depending on the   data and T.
  \end{theorem}

  A few remarks are in order.
\begin{remark} Our result can be regarded as a natural generalization of   Kazhikhov's theory (\cite{ka1})  for the constant viscosity and heat conductivity case to the degenerate and nonlinear ones.\end{remark}

  \begin{remark}Our  Theorem \ref{thm1.1} improves Hu-Ju's result  \cite{7} where they only treated the case that $\al=0,\b>0$ and   assumed that the initial data satisfy  \eqref{h1.8}   which is indeed   stronger than \eqref{a1.11}.  \end{remark}

  \begin{remark} Our  result still holds for compressible Navier-Stokes system ($\nb\equiv 0, \nw\equiv 0$) which generalized sightly  those due to \cite{hs1,28} where they only consider the case $\al=0,\b>0.$\end{remark}

We now make some comments on the analysis of this paper. To extend the local strong solutions whose existence can be obtained by   using the principle of compressed mappings (Lemma \ref{plq1}) to be global, compared with  \cite{7}, the key issue is to obtain the lower and upper bounds of both $v$ and $\te$ just under the conditions that the initial data satisfies \eqref{a1.11}.  Motivated by Kazhikhov \cite{ka1}, we first obtain a key representation of $v$    (see \eqref{2.6}). However,  if $\al>1,$ it seems difficult to obtain the lower bound of $v$ directly due to the nonlinearity of $\mu.$  To overcome this difficulty, we use the representation of $v,$ the energy-type inequality  \eqref{2.2}, and the Jensen's inequality to obtain a bound of $L^\infty(0,T;L^1)$-norm of $v^{-\al}$ (see \eqref{kjq2}) which plays an important role in bounding $v$ from below. Then,  after obtaining the estimates on the $L^2((0,1)\times(0,T))$-norm of both $\nb_{xx}$ and $\nw_{xx}$ (see \eqref{aas1}),  we multiply the momentum equation \eqref{1.2} by  $ u_{xx}$ and make full use of the structure of the energy equation \eqref{1.5}  to find that the $L^2((0,1)\times (0,T))$-norm of $u_{xx}$ can be bounded by the $L^2((0,1)\times (0,T))$-norm of $\te^{\b/2}\te_x$ (see \eqref{eee}) which indeed can be obtained by combining  the equation of $\te$ (see \eqref{2.4})  multiplied by $\te  $ and using the estimates obtained above (see \eqref{jjj}).  Once we get the bounds on the $L^2((0,1)\times (0,T))$-norm of both $u_{xx}$ and $u_t$ (see \eqref{aaas1}), the desired estimates on $\te_t$ and $\te_{xx}$ can be  obtained by standard arguments (see \eqref{eq1}). The details     will be carried out in the next section.

\section{ Proof of \thmref{thm1.1}}

 We first state   the following  existence and uniqueness of local solutions which can be obtained by using the Banach theorem and the contractivity of the operator defined by the linearization of the problem on a small time interval (c.f. \cite{10,13,tan}).

\begin{lemma} \la{plq1}Let \eqref{a1.10}-\eqref{a1.12} hold. Then there exists some $T>0$ such that  the initial-boundary-value problem \eqref{1.1}-\eqref{1.9} has a unique strong solution $(v,u,\te)$ satisfying\be
 \begin{cases}   v ,\,\theta \in L^\infty(0,T;H^1(0,1)),\quad u,\,\nb,\nw  \in L^\infty(0,T;H^1_0(0,1)),\\ v_t\in
  L^\infty(0,T;L^2(0,1))\cap L^2(0,T;H^1(0,1)), \\ u_t,\,\theta_t,\,\nb_t,\,\nw_t,\,u_{xx},\,\te_{xx},\,\nb_{xx},\,\nw_{xx} \in
  L^2((0,1)\times(0,T)),\end{cases} \ee   \end{lemma}

Theorem \ref{thm1.1} will be proved by extending the local solutions globally in time based on the global a priori estimates of solutions (see Lemma \ref{lemma30}--\ref{lemma70}) which will be obtained below.

Without loss of generality, we assume that $\lambda=\nu=\ti\mu_1= \ti\ka=R=c_v=1 , \ti\mu_2=\al,$ and that
 \bnn  \int_0^1 v_0dx=1.\enn

Then, we derive the following representation of $v$ which is essential in obtaining the  upper and lower bounds of $v$.
\begin{lemma} \la{lemma20} The  following expression of $v$ holds
\be \la{2.6}\ba v(x,t)=  B_0(x)  D (x,t)Y (t)  \left\{1 + \frac{1}{B_0(x)}\int_0^t\frac{(\te+  \frac{v}{2}|\nb|^2)(x,\tau)}{  D (x,\tau)Y (\tau)} d\tau \right\}, \ea\ee where
\be\ba\la{heq2.5}  B_0(x)= v_0 \exp\left\{- v_0^{-\al}- \int_0^1f_\al(v_0 )dx\right\}  ,\ea\ee
\be\ba\la{a2.7}D (x,t)=&\exp\left\{     v(x,t)^{-\al}+ \int_{0}^x \left(u(y,t)-u_0(y)\right)dy\right\}\\&\times \exp\left\{- \int_0^1v\int_0^xudydx+ \int_0^1v_0\int_0^xu_0dydx\right\},\ea\ee

\be\ba\la{2.8} Y (t)=\exp\left\{ \int_0^1f_\al(v  )dx- \int_0^t\int_0^1\left(u^2+ \frac{v}{2}|\nb|^2+\te\right) dxd\tau\right\} ,\ea\ee   with \be \la{d2.8}f_\al(s)=\begin{cases}\frac{\al}{1-\al}s^{1-\al} , &\mbox{ if } \alpha\not=1,\\ \ln s , &\mbox{ if } \alpha= 1 . \end{cases} \ee
\end{lemma}

 \pf  First, it follows from \eqref{1.2} that
\be\la{2.10}\ba u_t=\si_x,\ea\ee  where
\be\la{2.9}\ba  \si\triangleq \mu \frac{u_x}{v}-\frac{\te}{v}-\frac12|\nb|^2, \ea\ee satisfies \be \la{a2.9} \si=( \ln v-  v^{ -\al} )_t -\frac{\te}{v}-\frac12|\nb|^2,\ee due to \eqref{1.1}.
Integrating   \eqref{2.10}   over $(0,x)$ gives
\be\la{2.11}\ba \left(\int_0^xudy\right)_t=\si-\si(0,t) , \ea\ee
 which implies
\bnn\ba\la{2.0} v\si(0,t)=v\si-v\left(\int_0^xudy\right)_t . \ea\enn
 Then, integrating this in $x $ over $(0,1)$ and noticing that integrating \eqref{1.1} over $(0,1)\times (0,t)$ yields that   for any $t>0$
 \bnn\int_0^1 v(x,t)dx=  1,   \enn
 we obtain after using \eqref{1.9}    and \eqref{2.9}  that \be\ba\la{2.13}   \si(0,t) &= \int_0^1\left(\mu u_x-\te-\frac{v}{2}|\nb|^2 \right) dx-\left(\int_0^1v\int_0^xudydx\right)_t\\&\quad +\int_0^1u_x\int_0^xudy dx\\&= \left(  \int_0^1f_\al(v  )dx -\int_0^1v\int_0^xudydx\right)_t -\int_0^1\left(\te+ \frac{v}{2}|\nb|^2 + u^2 \right) dx.\ea\ee

 Finally, combining \eqref{2.11}, \eqref{a2.9}, and \eqref{2.13} yields
\bnn\la{2.12}\ba v(x,t) = B_0(x)  D (x,t)Y (t)  \exp\left\{ \int_0^t\left(\te+  \frac{v}{2}|\nb|^2\right)v^{-1} d\tau\right\}, \ea\enn
 with $B_0(x),$ $D (x,t),$ and $Y (t)$   as in  \eqref{heq2.5}-\eqref{2.8}  respectively. This in particular gives 
 \eqref{2.6} and finishes the proof of  Lemma  \ref{lemma20}.\thatsall

  With Lemma  \ref{lemma20} at hand,   we are in a position to prove
  lower bounds of both $v$ and $\te.$
\begin{lemma}\la{lemma30} It holds that for any $(x,t)\in[0,1]\times [0,T],$
\be\ba C^{-1}\le v(x,t) ,\quad
C^{-1}\le  \te(x,t),
\ea\la{2.15}\ee where (and in what follows)   $C $   denotes
some generic positive constant
 depending only on $T,\al, \b,\|(v_0 ,u_0,\theta_0 ,\nb_0,\nw_0)\|_{H^1(0,1)},
 \inf\limits_{x\in [0,1]}v_0(x),$ and $ \inf\limits_{x\in [0,1]}\theta_0(x).$
\end{lemma}

\pf First, using  \eqref{1.1}-\eqref{1.4}, we rewrite the energy equation \eqref{1.5}   as
\be\la{2.4}\theta_{t}+  \frac\theta v
u_{x}= \left(\frac{\theta^\b\theta_{x}}{v}\right)_{x}+ \frac{\mu u_{x}^{2}+|\nw_x|^2+|\nb_x|^2}{v}.\ee    Multiplying   \eqref{1.1},
   \eqref{1.2},       \eqref{1.3}, \eqref{1.4},   and \eqref{2.4}   by $ 1- {v}^{-1} ,
u, \nw, \nb,$ and  $  1- {\theta}^{-1} $  respectively,   adding them altogether and integrating the result over ${(0,1)\times(0,t)}$, we obtain the following energy-type inequality
 \be\ba\la{2.2} &\sup_{0\le t\le T}\int_0^1\left(  u^2+|\nw|^2+v|\nb|^2  + (v-\ln
v )+ (\theta-\ln \theta )\right)dx \\&\quad+\int_0^T V(s)ds \le C, \ea\ee
where
\bnn\ba V(t)\triangleq\int_0^1\left(\frac{\theta^\b\theta_{x}^2}{v\te^2}+\frac{\mu u_{x}^{2}}{v\te}
+\frac{|\nw_{x}|^{2}}{v\te}+\frac{|\nb_{x}|^{2}}{v\te}\right)(x,t)dx. \ea\enn  

 Next,    \eqref{2.2}   implies\bnn \left|\int_0^1 v\int_0^x udydx\right|  \le\int_0^1 v\left(\int_0^1 u^2dy\right)^{1/2}dx  \le C,
\enn which combined   with  \eqref{a2.7}  and \eqref{2.2} gives
\bn\la{2.18}\ba C^{-1} \le C^{-1}\exp\left\{     v(x,t)^{-\al}\right\} \le D(x,t)\le  C\exp\left\{     v(x,t)^{-\al}\right\} .\ea\en

Furthermore,  one deduces from   \eqref{2.2} that\be\la{d2.9}\ba  \left|\int_0^1 \ln vdx\right|+\int_0^1\left(u^2+ \frac{v}{2}|\nb|^2+\te\right) dx\le C,\ea\ee
which   yields that
\be\la{2.19}\ba C^{- 1 }\exp\left\{ \int_0^1f_\al(v  )dx\right\}\le  Y (t) \leq C\exp\left\{ \int_0^1f_\al(v  )dx\right\}  .\ea\ee
 Combining \eqref{2.6},  \eqref{2.18}, and  \eqref{2.19} yields that for any $(x,t)\in[0,1]\times[0,T],$
\be\la{2.20}\ba v(x,t )\geq   C^{- 1 }\exp\left\{  v(x,t)^{-\al}+ \int_0^1f_\al(v  )dx\right\},\ea\ee
which together with \eqref{d2.8} and \eqref{d2.9} leads to
\be\la{2.20k}\ba  \min_{(x,t)\in[0,1]\times[0,T]}v(x,t)\ge C^{-1},\ea\ee provided $\al\in [0,1].$
On the other hand, if $\al>1,$ integrating \eqref{2.20} in $x$ over $(0,1)$ and using  \eqref{d2.8}, \eqref{2.2}, and Jensen's inequality gives
\bnn \ba C&\ge  C^{- 1 } \exp\left\{  \int_0^1v(x,t)^{-\al}dx+\int_0^1f_\al(v  )dx\right\}\\&\ge  C^{- 1 } \exp\left\{ \frac12 \int_0^1v(x,t)^{-\al}dx-C\right\},\ea\enn
which in particular implies
\be\la{kjq2} \sup_{0\le t\le T}\int_0^1v(x,t)^{-\al}dx\le C.\ee
Hence, putting this into \eqref{2.20} shows \eqref{2.20k} still holds for $\al>1.$

 Finally, for $p>2,$ multiplying \eqref{2.4} by $  \te^{-p}$ gives
\bnn\ba\la{2.24k} & \frac{1}{p-1}\frac{d}{dt}\int_0^1\left( {\te}^{-1}\right)^{p-1}dx+\int_0^1\frac{\mu u_x^2}{v\te^p}dx \\&\le \int_0^1\frac{u_x}{v\te^{p-1}}dx\\&\leq
\frac{1}{2}\int_0^1\frac{\mu u_x^2}{v\te^p}dx+\frac{1}{2}\int_0^1\frac{1}{\mu v\te^{p-2}}dx\\&\leq\frac{1}{2}\int_0^1\frac{\mu u_x^2}{v\te^p}dx+C\left\|  \te^{-1} \right\|^{p-2}_{L^{p-1}} ,
 \ea\enn
where in the second inequality we have used  $\mu v \geq C^{-1}$. Combining this with  Gronwall's inequality yields that there exists some $C$ independent of $p$ such that
\bnn\ba \sup_{0\le t\le T}\left\|  \te^{-1}(\cdot,t) \right\|_{L^{p-1}}\leq C.\ea\enn  Letting $p\rightarrow\infty$  proves   the second inequality of \eqref{2.15}  and  finishes the proof of \lemref{lemma30}.\thatsall

\begin{lemma} \la{lemma40}  There exists a positive constant $C$ such that  for each $(x,t)\in [0,1]\times[0,T],$
\be\ba\la{aee} C^{-1}\leq v(x,t)\leq   C.\ea\ee
\end{lemma}

\pf
First, for $0<\alpha< 1 $ and $0<\ve<1,$  integrating \eqref{2.4} multiplied by $\te^{-\al}$ over $ (0,1)\times(0,T)$ yields
\be\ba\la{b}&\int_0^T\int_0^1\frac{\al\te^\b\te_x^2}{v\te^{\al+1}}dxdt+\int_0^T\int_0^1\frac{\mu u_x^2 +|\nw_x|^2+|\nb_x|^2}{v\te^\al}dxdt \\& =  \frac{1}{1-\al}\int_0^1\left(\te^{1-\al}-\te_0^{1-\al}\right)dx+\int_0^T\int_0^1\frac{\te^{1-\al}u_x}{v} dx dt \\&  \leq C(\al)+\frac{1}{2}\int_0^T\int_0^1\frac{\mu u_x^2}{v\te^\al}dxdt+C\int_0^T\int_0^1 \te^{2-\al} dxdt \\&\leq C(\al)+\frac{1}{2}\int_0^T\int_0^1\frac{\mu u_x^2}{v\te^\al}dxdt
+C\int_0^T\max_{x\in[0,1]}\te^{1-\al}\int_0^1\te dxdt \\&\leq C(\al,\ve)+\frac{1}{2}\int_0^T\int_0^1\frac{\mu u_x^2}{v\te^\al}dxdt+\ve \int_0^T\max_{x\in[0,1]}\te dt,\ea\ee where in the first inequality we have used \eqref{2.2} and \eqref{2.15}.

Next,  for $ \al=\min\{1,\b\}/2,$ using   \eqref{2.2}, we get
\bnn\ba\la{f}\int_0^T\max_{x\in[0,1]}\te dt&\leq C+C\int_0^T\int_0^1|\te_x|dxdt \\&\leq C+C\int_0^T\int_0^1\frac{ \te^\b\te_x^2}{v\te^{1+\al}}dxdt+C\int_0^T\int_0^1\frac{v\te^{1+\al}}{\te^\b}dxdt \\&\leq C +C\int_0^T\int_0^1\frac{ \te^\b\te_x^2}{v\te^{1+\al}}dxdt+\frac{1}{2}\int_0^T\max_{x\in[0,1]}\te dt,\ea\enn which together with \eqref{b} yields that
\bn\la{g}\ba\int_0^T\max_{x\in[0,1]}\te dt\leq C  ,\ea\en and then that
for $0<\alpha< 1, $
\be\la{aa}\ba\int_0^T\int_0^1\frac{ \te^\b\te_x^2}{v\te^{\alpha+1}}dxdt \leq C(\al).\ea\ee

Finally, it follows from \eqref{2.19}, \eqref{d2.8},    \eqref{d2.9}, \eqref{2.2},  and   \eqref{2.20k} that \bnn C^{-1}\le Y(t)\le C,\enn which together with   \eqref{2.6},    \eqref{2.18},  \eqref{2.20k},  and   \eqref{g} yields
\be \la{j1}v(x,t)\le C+C\int_0^t \max_{x\in[0,1]}|\nb|^2(x,t)\max_{x\in[0,1]} v(x,t)dt.\ee
Using
\eqref{1.9},  \eqref{2.2}, and  \eqref{g}, we have
\be \la{44}\ba  \int_0^T\max_{x\in[0,1]}|\nb|^2(x,t)dt&\leq C\int_0^T\int_0^1|\nb\cdot\nb_x| dxdt \\&\leq C\int_0^T\int_0^1\frac{|\nb_x|^2}{v\te}dxdt+C\int_0^T\int_0^1v\te|\nb|^2dxdt \\&\leq C+C\int_0^T\max_{x\in[0,1]}\te dt \\&\leq C,\ea\ee
   which combined with  \eqref{j1}  and Gronwall's inequality  gives
\bn\la{11}\ba\max_{(x,t)\in[0,1]\times[0,T]}v\leq C .  \ea\en The   proof of \lemref{lemma40} is finished. \thatsall
 \begin{lemma}\la{lemma60}
There is a positive constant C such that
\begin{align}\la{aaaa}
&\sup_{0\le t\le T} \int_0^1 v_x^2dx \leq C.\end{align}
\end{lemma}

\pf First, we rewrite the momentum equation \eqref{1.2} as
\bnn\la{bbbb}\ba \left(u-\frac{\mu v_x}{v}\right)_t=-\left(\frac{\te}{v}+\frac{1}{2}|\nb|^2\right)_x.   \ea\enn
Multiplying the above equation by $u-\frac{\mu v_x}{v}$ and integrating the resultant equality  yields that for any $t\in(0,T)$
\be\ba\la{cccc}&\frac{1}{2}\int_0^1\left(u-\frac{\mu v_x}{v}\right)^2 dx -\frac{1}{2}\int_0^1\left(u-\frac{\mu v_x}{v}\right)(x,0)dx\\&= \int_0^t \int_0^1\left( \frac{\te v_x}{v^2}-\frac{\te_x}{v}-\nb\cdot\nb_x\right)\left(u-\frac{\mu v_x}{v}\right)dxdt \\
&=  -\int_0^t\int_0^1\frac{\mu \te v_x^2}{v^3}dxdt+\int_0^t\int_0^1\frac{\te u v_x}{v^2}dxdt \\&\quad-\int_0^t\int_0^1\frac{\te_x}{v} \left(u-\frac{\mu v_x}{v}\right)dxdt-\int_0^t\int_0^1\nb\cdot\nb_x  \left(u-\frac{\mu v_x}{v}\right)dxdt \\
&=  -\int_0^t\int_0^1\frac{\mu \te v_x^2}{v^3}dxdt+\sum_{i=1}^3I_i.\ea\ee
Each $I_i (i=1,2,3)$ can be estimated as follows:

First, Cauchy's inequality gives
\be\ba\la{dddd}|I_1|&\leq\frac{1}{2}\int_0^t\int_0^1\frac{\mu \te v_x^2}{v^3}dxdt+\frac{1}{2}\int_0^T\int_0^1\frac{u^2\te}{\mu v}dxdt \\&\leq\frac{1}{2}\int_0^t\int_0^1\frac{\mu \te v_x^2}{v^3}dxdt+C\int_0^T\max_{x\in[0,1]}\te dt \\&\leq C+\frac{1}{2}\int_0^t\int_0^1\frac{\mu \te v_x^2}{v^3}dxdt,\ea\ee
where  we have used \eqref{aee}, \eqref{2.2}, and \eqref{g}.

Next, using \eqref{2.15}, \eqref{2.2} and \eqref{aee}, we have
\be\ba\la{eeee}|I_2| &\leq\frac{1}{2}\int_0^T\int_0^1\frac{\te^\b\te_x^2}{v\te^2}dxdt+\frac{1}{2}\int_0^t\int_0^1 \frac{\te^2}{v\te^\b}\left(u-\frac{\mu v_x}{v}\right)^2dxdt \\
&\leq C+C\int_0^t\max_{x\in[0,1]}\te^{2 }\int_0^1\left(u-\frac{\mu v_x}{v}\right)^2dxdt.\ea\ee
Moreover, it follows from   \eqref{2.2}, \eqref{2.15}, \eqref{aee}, and \eqref{aa} that, for any $\ve>0,$
\bnn\ba \int_0^T\max_{x\in[0,1]}\te^2 dt&\le  C\int_0^T\max_{x\in[0,1]}\left|\te^2-\int_0^1\te^2dx\right| dt+C\int_0^T \max_{x\in[0,1]}\te dt\\ &\leq C+C\int_0^T\int_0^1\te|\te_x| dxdt \\&\leq C+C(\ve) \int_0^T\int_0^1\frac{ \te_x^2}{v\te} dx dt+ \varepsilon\int_0^T\int_0^1v\te^3dxdt \\&\leq C(\ve)+C\varepsilon\int_0^T\max_{x\in[0,1]}\te^2dt,\ea\enn  which gives
\bn\la{22}\ba  \int_0^T\max_{x\in[0,1]}\te^2 dt\leq C .\ea\en

 Finally, integrating \eqref{2.4} over $(0,1)\times (0,T),$ we have by \eqref{aee}
\bnn \ba&\int_0^T\int_0^1\frac{\mu u_x^2+|\nw_x|^2+|\nb_x|^2}{v}dxdt\\ &= \int_0^1\te dx-\int_0^1\te_0dx +\int_0^T\int_0^1\frac{\te}{v}u_xdx\\&\le C+\frac{1}{2}\int_0^T\int_0^1\frac{\mu u_x^2 }{v}dxdt +C \int_0^T\max_{x\in[0,1]}\te^2 dt , \ea\enn
which together with \eqref{aee} and  \eqref{22} gives
\be\la{h}\ba&\int_0^T\int_0^1(u_x^2+|\nw_x|^2+|\nb_x|^2)dxdt  \le C . \ea\ee Combining this with Cauchy's inequality leads to
\be\ba\la{gggg}|I_3|&\leq C\int_0^t\int_0^1\left(|\nb_x|^2+|\nb|^2 \left(u-\frac{\mu v_x}{v}\right)^2\right)dxdt \\&\leq C+C\int_0^t\max_{x\in[0,1]}|\nb|^2\int_0^1\left(u-\frac{\mu v_x}{v}\right)^2dxdt.\ea\ee
Putting \eqref{dddd}, \eqref{eeee}, and \eqref{gggg} into \eqref{cccc}, we  obtain after using  Gronwall's inequality,  \eqref{44}, and \eqref{22}  that
\bnn\la{nnnn}\ba \sup_{0\le t\le T}\int_0^1\left(u-\frac{\mu v_x}{v}\right)^2dx+\int_0^T\int_0^1\frac{\te v_x^2}{v^3}dxdt\leq C,\ea\enn
which together with \eqref{2.2}  gives \eqref{aaaa} and finishes the proof of Lemma \ref{lemma60}. \thatsall

\begin{lemma}\la{lemmy0}
There is a positive constant C such that
\be \la{aas1}\ba
&\sup_{0\le t\le T} \int_0^1\left( |\nb_x|^2+|\nw_x|^2\right)dx \\&  +\int_0^T\int_0^1\left(  |\nb_t|^2 +|\nb_{xx}|^2 +|\nw_t|^2+|\nw_{xx}|^2\right)dx dt\leq C.\ea\ee
\end{lemma}

\pf First, multiplying \eqref{1.3} by $\nw_{xx}$ and integrating the resulting equality over $(0,1)\times(0,T)$, we obtain after using \eqref{1.9}, \eqref{h}, \eqref{aaaa}, and Cauchy's inequality that
\be\ba\la{bbb}&\frac{1}{2}\int_0^1|\nw_x|^2dx+\int_0^T\int_0^1\frac{|\nw_{xx}|^2}{v}dxdt   \\& \leq C+\frac{1}{2}\int_0^T\int_0^1\frac{|\nw_{xx}|^2}{v}dxdt +C\int_0^T\int_0^1\left(|\nb_x|^2+|\nw_x|^2v_x^2\right)dxdt \\&\leq C+\frac{1}{2}\int_0^T\int_0^1\frac{|\nw_{xx}|^2}{v}dxdt +C\int_0^T\max_{x\in[0,1]}|\nw_x|^2dt.\ea\ee
Direct computation shows  for any $\ve>0,$
\be\ba\la{ccc}\int_0^T\max_{x\in[0,1]}|\nw_x|^2dt&\leq C(\ve)\int_0^T\int_0^1|\nw_x|^2dxdt+\varepsilon\int_0^T\int_0^1\frac{|\nw_{xx}|^2}{v}dxdt \\&\leq C(\ve)+\varepsilon\int_0^T\int_0^1\frac{|\nw_{xx}|^2}{v}dxdt,\ea\ee
which combined with \eqref{bbb}  leads to
\bn\ba \la{ddd}\sup_{0\le t\le T}\int_0^1|\nw_x|^2dx+\int_0^T\int_0^1|\nw_{xx}|^2dxdt\leq C.\ea\en

Then, we rewrite  \eqref{1.3} as
\bnn\la{ppp}\ba \nw_t=\frac{\nw_{xx}}{v}-\frac{\nw_xv_x}{v^2}+\nb_x, \ea\enn
 which together with \eqref{aee}, \eqref{ddd}, \eqref{aaaa}, \eqref{h}, and \eqref{ccc} gives
\be\ba\la{ooo}\int_0^T\int_0^1|\nw_t|^2dxdt &\leq
C\int_0^T\int_0^1 \left(|\nb_x|^2+|\nw_{xx}|^2 + v_x^2|\nw_x|^2\right)dxdt\\&\leq C\int_0^T\max_{x\in[0,1]}|\nw_x|^2dt\\&\leq C .\ea\ee

Next, multiplying \eqref{1.4} by $\nb_{xx}$ and integrating the result over $(0.1)\times(0,T)$, we deduce from \eqref{aaaa},  \eqref{h}, \eqref{aee}, \eqref{2.2} and Cauchy's inequality that
\bnn\ba &\frac{1}{2}\int_0^1|\nb_x|^2dx+\int_0^T\int_0^1\frac{|\nb_{xx}|^2}{v^2}dxdt\\ &\leq C+\frac{1}{2}\int_0^T\int_0^1\frac{|\nb_{xx}|^2}{v^2}dxdt +C\int_0^T\int_0^1\left(|\nb_x|^2v_x^2+u_x^2|\nb|^2+|\nw_x|^2\right)dxdt \\&\leq
C+\frac{1}{2}\int_0^T\int_0^1\frac{|\nb_{xx}|^2}{v^2}dxdt+C\int_0^T\max_{x\in[0,1]}|\nb_x|^2dt +\max_{(x,t)\in[0,1]\times[0,T]}|\nb|^2 \\&\leq C+\frac{3}{4}\int_0^T\int_0^1\frac{|\nb_{xx}|^2}{v^2}dxdt+C \int_0^T\int_0^1|\nb_x|^2dxdt \\&\quad  +C\sup_{0\le t\le T}\int_0^1|\nb|^2dx +\frac{1}{4}\sup_{0\le t\le T}\int_0^1|\nb_x|^2dx \\&\leq C +\frac{3}{4}\int_0^T\int_0^1\frac{|\nb_{xx}|^2}{v^2}dxdt+\frac{1}{4}\sup_{0\le t\le T}\int_0^1|\nb_x|^2dx,\ea\enn which implies
\bn\ba\la{zzz}\sup_{0\le t\le T}\int_0^1|\nb_x|^2dx+\int_0^T\int_0^1|\nb_{xx}|^2dxdt\leq C.\ea\en
Hence, \be\la{xxx} \max_{(x,t)\in [0,1]\times [0,T]}|\nb|^2\le C+C \sup_{0\le t\le T}\int_0^1|\nb_x|^2dx\le C.\ee

Finally, we rewrite \eqref{1.4} as
\bnn\la{kkk}\ba\nb_t=\frac{\nw_x}{v}+\frac{\nb_{xx}}{v^2}-\frac{\nb_x v_x}{v^3}-\frac{\nb u_x}{v},\ea\enn
which together with  \eqref{zzz}, \eqref{aaaa}, \eqref{h}, and  \eqref{xxx}  gives
\bnn\ba\la{rrr}\int_0^T\int_0^1|\nb_t|^2dxdt&\leq C\int_0^T\int_0^1\left(|\nb_{xx}|^2+|\nb_x|^2v_x^2+|\nw_x|^2+|\nb|^2u_x^2\right)dxdt\\&\leq C+C\int_0^T\left(\max_{x\in[0,1]}|\nb_x|^2+ \int_0^1u_x^2dx\right)dt\\&\leq C+C\int_0^T\int_0^1\left(|\nb_x|^2+|\nb_{xx}|^2\right)dxdt\\&\leq C.\ea\enn
 Combining this,  \eqref{ddd}, \eqref{ooo},  and  \eqref{zzz}  gives \eqref{aas1} and finishes the proof of Lemma \ref{lemmy0}. \thatsall

\begin{lemma}\la{lemmy11}
There is a positive constant C such that
\be \la{aaas1}\ba
&\sup_{0\le t\le T} \int_0^1 u_x^2 dx   +\int_0^T\int_0^1\left( u_t^2 +u_{xx}^2 \right)dx dt\leq C.\ea\ee
\end{lemma}

\pf
First, multiplying \eqref{1.2} by $u_{xx}$ and integrating the result over $(0,1)\times(0,T)$,  we have
 \be\ba\la{eee} &\frac{1}{2}\int_0^1u_x^2dx+\int_0^T\int_0^1\frac{\mu u_{xx}^2}{v}dxdt  \\&\leq C+\frac{1}{2}\int_0^T\int_0^1\frac{\mu u_{xx}^2}{v}dxdt+C\int_0^T\int_0^1\left(\te_x^2+\te^2v_x^2+|\nb|^2|\nb_x|^2+u_x^2v_x^2 \right)dxdt \\&\leq C+\frac{1}{2}\int_0^T\int_0^1\frac{\mu u_{xx}^2}{v}dxdt+C\int_0^T\int_0^1\te_x^2dxdt +C\int_0^T\max_{x\in[0,1]}\te^2\int_0^1v_x^2dxdt \\&\quad
+C\max_{(x,t)\in[0,1]\times[0,T]}|\nb|^2\int_0^T\int_0^1|\nb_x|^2dxdt
+C\int_0^T\max_{x\in[0,1]}u_x^2\int_0^1v^2_xdxdt  \\&\leq C+ \frac{3}{4} \int_0^T\int_0^1\frac{\mu u_{xx}^2}{v}dxdt
+C_1\int_0^T\int_0^1\frac{\te^\b\te_x^2}{v}dxdt, \ea\ee
where in the last inequality we have used \eqref{22}, \eqref{aaaa}, \eqref{xxx}, \eqref{h}  and the following  inequality,
\be\ba\la{fff}\int_0^T\max_{x\in[0,1]}u_x^2dt&\leq
C(\ve)\int_0^T\int_0^1u_x^2dxdt+\ve\int_0^T\int_0^1\frac{\mu u_{xx}^2}{v}dxdt\\&\leq C(\ve)+\ve\int_0^T\int_0^1\frac{\mu u_{xx}^2}{v}dxdt,\ea\ee for any $\ve>0.$

Then, multiplying \eqref{2.4} by $\te$ and integrating the result over $(0,1)\times(0,T)$ yields
\be\ba\la{iii}&\frac{1}{2}\int_0^1\te^2dx+\int_0^T\int_0^1\frac{\te^\b\te_x^2}{v}dxdt\\&\leq C+C \int_0^T\int_0^1\te^2|u_x|dx+C\int_0^T\int_0^1\left( u_x^2 + |\nw_x|^2 + |\nb_x^2| \right)\te dxdt\\&\leq C+C\int_0^T\int_0^1\te u_x^2dxdt+C\int_0^T\int_0^1\te^3dxdt+\int_0^T\max_{x\in[0,1]}\te dt\\&\leq C+C\int_0^T \max_{[0,1]}u_x^2 dt+C\int_0^T\max_{x\in[0,1]}\te^2dt \\&\leq C(\varepsilon)+C\ve\int_0^T\int_0^1\frac{\mu u_{xx}^2}{v}dxdt, \ea\ee
 where we have used \eqref{aas1}, \eqref{fff} and \eqref{22}. Adding \eqref{iii} multiplied by $ C_1+1$  to \eqref{eee} and choosing $\varepsilon$  sufficiently small, we obtain  that
\bn\la{jjj}\ba \sup_{0\le t\le T}\int_0^1(\te^2 + u_x^2)dx+\int_0^T\int_0^1 \te^\b\te_{x}^2 dxdt+\int_0^T\int_0^1 u_{xx}^2   dxdt\leq C. \ea\en

 Finally, we rewrite \eqref{1.2} as
\bnn\ba\la{eq20} u_t=\frac{\mu u_{xx}}{v}-\left(\frac{\mu }{v }\right)'_ v u_xv_x-\frac{\te_x}{v}+\frac{\te v_x}{v^2}-\nb\cdot\nb_x ,\ea\enn
which together with \eqref{aaaa}, \eqref{xxx}, \eqref{22}, \eqref{fff}, \eqref{jjj},  and \eqref{aas1} leads to
\bnn\ba
\int_0^T\int_0^1u_t^2dxdt&\leq C\int_0^T\int_0^1\left(u_{xx}^2+u_x^2v_x^2+\te_x^2+\te^2v_x^2+|\nb|^2|\nb_x|^2\right)dxdt\\&\leq C.
\ea\enn
Combining this and   \eqref{jjj}  immediately gives \eqref{aaas1} and  completes the proof of \lemref{lemmy11}.\thatsall

\begin{lemma}\la{lemma70}There exists a positive constant $C$ such that \bn\ba\la{eq1} \sup_{0\le t\le T}\xix  \te_x^2dx+\int_0^T\xix \left( \te_t^2+\te_{xx}^2\right)dxdt\le C .  \ea\en
\end{lemma}

\pf
  First, noticing that integration by parts leads to
\bnn\ba\la{eq3}  \xix \te^{\b}\te_t\left(\frac{\te^\b\theta_{x}}{v}\right)_{x}dx & =-\xix \frac{\te^\b\theta_{x}}{v}\left(\te^{\b}\te_t\right)_{x}dx\\& =-\xix \frac{\te^\b\theta_{x}}{v}\left(\te^{\b}\te_x\right)_{t}dx \\&=-\frac{1}{2}
\xix \frac{\left((\te^\b\theta_{x})^2\right)_t}{v}dx \\&=-\frac{1}{2} \left(\xix \frac{(\te^\b\theta_{x})^2}{v}dx\right)_t
-\frac{1}{2}\xix \frac{(\te^\b\theta_{x})^2u_x}{v^2}dx  ,\ea\enn
multiplying \eqref{2.4} by
$ \te^{\b}\te_t$ and integrating the resultant equality over (0,1), we have
\be\ba\la{eq4} & \xix  \te^\b\te_t^2dx+ \frac{1}{2} \left(\xix \frac{(\te^\b\theta_{x})^2}{v}dx\right)_t\\&=- \frac{1}{2}\xix \frac{(\te^\b\theta_{x})^2u_x}{v^2}dx+\xix  \frac{ \te^{\b }\theta_t\left(-\te u_x+\mu u_x^2+|\nw_x|^2+|\nb_x|^2\right)}{v }dx\\&\le C\max_{x\in[0,1]} (|u_x|\te^{\b/2})\xix \te^{3\b/2}\te_x^2dx +\frac{1}{2}\xix  \te^\b\te_t^2dx+C\xix  \te^{\b+2}u_x^2dx\\&\quad+C\xix \te^\b\left(u_x^4+|\nw_x|^4+|\nb_x|^4\right)dx
\\&\le  C \xix  \te^{2\b}\te_x^2dx\xix  \te^\b\te_x^2dx+\frac{1}{2}\xix  \te^\b\te_t^2dx\\&\quad+C \max_{x\in[0,1]}(  \te^{2\b+2}+  u_x^4 +|\nw_x|^4+|\nb_x|^4 )+C ,\ea\ee due to
\eqref{2.15},  \eqref{aas1}  and \eqref{jjj}.

Combining \eqref{aaas1} with H\"{o}lder's inequality gives
 \be\ba\la{eq5}\int_0^T\max_{x\in[0,1]}u_x^4dt&\leq C\int_0^T\int_0^1u_x^4dxdt+C\int_0^T\int_0^1|u_x^3u_{xx}|dxdt\\&\leq C\int_0^T\max_{x\in[0,1]}u_x^2\int_0^1u_x^2dxdt\\&\quad+C\int_0^T\max_{x\in[0,1]}u_x^2
 \left(\int_0^1u_x^2dx\right)^{\frac{1}{2}}\left(\int_0^1u_{xx}^2dx\right)^\frac{1}{2}dt\\&\leq C\int_0^T\int_0^1\left(u_x^2 + u_{xx}^2\right)dxdt +\frac12 \int_0^T \max_{x\in[0,1]}u_x^4dt
 \\&\leq C +\frac12 \int_0^T \max_{x\in[0,1]}u_x^4dt .\ea\ee
  Using \eqref{aas1} and applying similar arguments to $\nb$ and $\nw$ implies
 \bn\la{eq6}\ba\int_0^T\max_{x\in[0,1]}|\nb_x|^4dt\leq C,\quad \int_0^T\max_{x\in[0,1]}|\nw_x|^4dt\leq C. \ea\en
Noticing that \be\la{eq9}\max_{x\in[0,1]} \te^{2\beta+2} \le C+C\xix (\te^\b\te_x)^2dx,\ee we then deduce from \eqref{jjj}, \eqref{eq4}--\eqref{eq6}, and the Gronwall inequality that
\be\ba\la{eq8} \sup_{0 \le t\le T}\xix  \left(\te^\b\theta_{x}\right)^2 dx+\int_0^T\xix  \te^\b\te_t^2dxdt\le C, \ea\ee
which together with \eqref{eq9}  shows
\be\la{eq10}\max_{(x,t)\in[0,1]\times[0,T]}\te(x,t)\le C.\ee
Thus, both \eqref{eq8} and \eqref{2.15} lead   to
\be\ba\la{eq11}\sup_{0 \le t\le T}\xix  \theta_{x}^2 dx+\int_0^T\xix  \te_t^2dxdt\le C. \ea\ee

Finally, it follows from \eqref{2.4} that
\bnn\ba \frac{\te^\b\te_{xx}}{v}= -\frac{\b \te^{\b-1}\te_x^2}{v}+\frac{\te^\b\te_x v_x}{v^2}- \frac{\mu u_x^2+|\nb_x|^2+|\nw_x|^2}{v}+ \frac{ \te u_x}{v}+\te_t,\ea\enn
which together with \eqref{2.15}, \eqref{h}, \eqref{aaaa}, \eqref{eq5}, \eqref{eq6}, \eqref{eq10}, and \eqref{eq11} yields
\bnn\ba\la{eq12}\int_0^T\xix \te_{xx}^2dxdt  & \le C\int_0^T\xix \left(\te_x^4+\te_x^2v_x^2+u_x^4+|\nb_x|^4+|\nw_x|^4+u_x^2+\te_t^2\right)dxdt\\ & \le C+ C\int_0^T\max_{x\in[0,1]}\te_x^2 dt\\ & \le C+   \frac12\int_0^T  \xix \te_{xx}^2 dx   dt.\ea\enn
 Combining this with \eqref{eq11}  shows \eqref{eq1}
and  finishes the proof of \lemref{lemma70}.\thatsall


 \end{document}